%&amstex     
\input amstex\documentstyle {amsppt}  
\pagewidth{12.5 cm}\pageheight{19 cm}\magnification\magstep1
\topmatter
\title A $q$-analogue of an identity of N.Wallach\endtitle
\author G. Lusztig\endauthor
\address Department of Mathematics, M.I.T., Cambridge, MA 02139\endaddress
\thanks Supported in part by the National Science Foundation\endthanks
\endtopmatter   
\document

\define\lb{\linebreak}

\define\m{\mapsto}
\define\do{\dots}

\define\sub{\subset}

\define\T{\times}

\define\nl{\newline}

\define\ot{\otimes}

\define\ph{\phi}

\redefine\t{\tau}

\define\qq{\bold q}

\redefine\tt{\bold t}

\define\CC{\bold C}

\define\FF{\bold F}

\define\QQ{\bold Q}

\define\ZZ{\bold Z}

\define\cf{\Cal F}

\define\ch{\Cal H}

\define\GW{GW}
\define\DFP{DFP}
\define\WA{W}
\define\PH{Ph}
\subhead 1\endsubhead
Let $n$ be an integer $\ge 2$. For $i\in[1,n-1]$ let $s_i$ be the transposition 
$(i,i+1)$ in the group $S_n$ of permutations of $[1,n]$. (Given two integers $a,b$ we 
denote by $[a,b]$ the set of all integers $c$ such that $a\le c\le b$.) Consider the 
following element of $\CC[S_n]$ (the group algebra of $S_n$):
$$\tt=s_1s_2\do s_{n-1}+s_2s_3\do s_{n-1}+\do+s_{n-2}s_{n-1}+s_{n-1}+1.$$
(The sum of an $n$-cycle,an $n-1$-cycle,..., a $2$-cycle and the identity.) Wallach
\cite{\WA} proved the remarkable identity 
$$\tt\prod_{k\in[1,n],k\ne n-1}(\tt-k)=0\tag a$$
in $\CC[S_n]$ and used it to establish a vanishing result for some Lie algebra 
cohomologies. In particular, left multiplication by $\tt$ in $\CC[S_n]$ has
eigenvalues in $\{0,1,2,\do,n-2,n\}$. 
A closely related result appeared later in connection with a problem 
concerning shuffling of cards in Diaconis, Fill and Pitman \cite{\DFP} and also in 
Phatarfod \cite{\PH}. 

Let $\qq$ be an indeterminate. Let $H$ be the $\ZZ[\qq]$-algebra with generators \lb
$T_1,T_2,\do,T_{n-1}$ and relations $(T_i+1)(T_i-\qq)=0$ for $i\in[1,n-1]$, 
$T_iT_{i+1}T_i=T_{i+1}T_iT_{i+1}$ for $i\in[1,n-2]$, $T_iT_j=T_jT_i$ for $i\ne j$ in 
$[1,n-1]$, a Hecke algebra of type $A_{n-1}$. Set
$$\t=T_1T_2\do T_{n-1}+T_2T_3\do T_{n-1}+\do+T_{n-2}T_{n_1}+T_{n-1}+1\in H.$$
Under the specialization $\qq=1$, $\t$ becomes the element $\tt$ of $\CC[S_n]$. Our 
main result is the following $q$-analogue of (a):

\proclaim{Proposition 2} The following equality in $H$ holds:
$$\t\prod_{k\in[1,n],k\ne n-1}(\t-1-\qq-\qq^2-\do-\qq^{k-1})=0.$$
\endproclaim   
The proof will be given in Section 4. The proof of the Proposition is a generalization 
of the proof of 1(a) given in \cite{\GW}. However, there is a new difficulty due to the
fact that the product of two standard basis elements of $H$ is not a standard basis 
element (as for $S_n$) but a complicated linear combination of basis elements. To 
overcome this difficulty we will work in a model of $H$ as a space of functions on a 
product of two flag manifolds over a finite field. 

Let $V$ be a vector space of dimension $n$ over a finite field $\FF_q$ of cardinal $q$.
Let $\cf$ be the set of complete flags 
$$V_*=(V_0\sub V_1\sub V_2\sub\do\sub V_n)$$ 
in $V$ where $V_k$ is a subspace of $V$ of dimension $k$ for $k\in[0,n]$. Now $GL(V)$ 
acts on $\cf$ by 
$$g:V_*\m gV_*=(gV_0\sub gV_1\sub gV_2\sub\do\sub gV_n)$$
and on $\cf\T\cf$ by $g:(V_*,V'_*)\m(gV_*,gV'_*)$. Let $\ch$ be the $\CC$-vector space 
of all functions $f:\cf\T\cf@>>>\CC$ that are constant on the orbits of $GL(V)$. This 
is an associative algebra with multiplication 
$$f,f'\m f*f',\quad(f*f')(W_*,V'_*)=\sum_{V'_*\in\cf}f(W_*,V'_*)f'(V'_*,V_*).$$
Define $f_1\in\ch$ by

$f_1(W_*,V'_*)=1$ if there exists $g\in[1,n]$ (necessarily unique) with $W_r=V'_r$ for 
$r\in[1,g-1]$, $V'_r\ne W_r\sub V'_{r+1}$ for $r\in[g,n-1]$;

$f_1(W_*,V'_*)=0$, otherwise.
\nl
For $t\in[0,n]$ and any sequence $1\le i_1<i_2<\do<i_{n-t}\le n$ let 
$X_t^{i_1,i_2,\do,i_{n-t}}$ be the set of all pairs $(V'_*,V_*)\in\cf\T\cf$ such that
$V'_r\sub V_{i_r},V'_r\not\sub V_{i_r-1}$ for $r\in[1,n-t]$. For $t\in[0,n]$ let 
$X_t=\cup X_t^{i_1,i_2,\do,i_{n-t}}\sub\cf\T\cf$ where the union is taken over all 
sequences $1\le i_1<i_2<\do<i_{n-t}\le n$. Clearly, this union is disjoint and
$X_0\sub X_1\sub X_2\sub\do\sub X_n=\cf\T\cf$. Also, $X_0$ is the diagonal in 
$\cf\T\cf$. Define $f_t\in\ch$ by

$f_t(V'_*,V_*)=1$ if $(V'_*,V_*)\in X_t$, $f_t(V'_*,V_*)=0$, otherwise.  
\nl
For $t=1$ this agrees with the earlier definition of $f_1$. Note that $f_0$ is the unit
element of $\ch$. The following result is a $q$-analogue of a result in \cite{\DFP}.

\proclaim{Lemma 3}For $t\in[1,n-1]$ we have 
$f_1*f_t=(1+q+q^2+\do+q^{t-1})f_t+q^tf_{t+1}$.
\endproclaim
Let $f=f_1*f_t$. From the definitions we have $f=\sum_{g=1}^n\ph_g$ where $\ph_g\in\ch$
is defined as follows: for $(W_*,V_*)\in\cf\T\cf$, $\ph_g(W_*,V_*)$ is the number of 
$V'_*\in\cf$ such that 

$V'_r=W_r$ for $r\in[1,g-1]$, $V'_r\ne W_r\sub V'_{r+1}$ for $r\in[g,n-1]$ and there 
exists $1\le i_1<i_2<\do<i_{n-t}\le n$ with $V'_r\sub V_{i_r},V'_r\not\sub V_{i_r-1}$ 
for $r\in[1,n-t]$. 
\nl
Here $V'_r$ is uniquely determined for $r\in[1,g-1]$ (we have $V'_r=W_r$) while for
$r\in[g+1,n-1]$, $V'_r$ is equal to $V'_g+W_{r-1}$ (this follows by induction from
$V'_r=V'_{r-1}+W_{r-1}$ which holds since $V'_{r-1},W_{r-1}$ must be distinct
hyperplanes of $V'_r$). Hence $\ph_g(W_*,V_*)$ is the cardinal of the set
$Y_g$ consisting of all $g$-dimensional subspaces $V'_g$ of $V$ such that 

$W_{g-1}\sub V'_g$,

$V'_g+W_{r-1}\ne W_r$ for $r\in[g,n-1]$ (or equivalently $V'_g\not\sub W_{n-1}$),
\nl
and there exists $1\le i_1<i_2<\do<i_{n-t}\le n$ (necessarily unique) with 

$W_r\sub V_{i_r},W_r\not\sub V_{i_r-1}$ if $r\in[1,n-t]\cap[1,g-1]$,

$V'_g\sub V_{i_g},V'_g\not\sub V_{i_g-1}$ if $g\in[1,n-t]$,

$V'_g+W_{r-1}\sub V_{i_r},V'_g+W_{r-1}\not\sub V_{i_r-1}$ if 
$r\in[1,n-t]\cap[g+1,n-1]$.
\nl
Assume first that $g\in[1,n-t]$. If a $V'_g\in Y_g$ exists 
and if $1\le i_1<i_2<\do<i_{n-t}\le n$ is as above then,
setting $j_r=i_r$ for $r\in[1,g-1]$ and $j_r=i_{r+1}$ for $r\in[g,n-t-1]$, 
we have $1\le j_1<j_2<\do<j_{n-t-1}\le n$ and 

(a) $W_r\sub V_{j_r},W_r\not\sub V_{j_r-1}$ for $r\in[1,n-t-1]$. 
\nl
(For $r\in[1,g-1]$ this is clear. Assume now that $r\in[g,n-t-1]$. Since 
$V'_g+W_r\sub V_{j_r}$, we have $W_r\sub V_{j_r}$. If $W_r\sub V_{j_r-1}$ then, since 
$V'_g\sub V_{i_g}\sub V_{j_r-1}$ and $j_r=i_{r+1}$, we would have 
$V'_g+W_r\sub V_{i_{r+1}-1}$, contradiction.) We see that $\ph_g(W_*,V_*)=0$ if 
$(W_*,V_*)\notin X_{t+1}$. We now assume that 
$(W_*,V_*)\in X_{t+1}$. Let $1\le j_1<j_2<\do<j_{n-t-1}\le n$ be such that (a) holds. 
Then $\ph_g(W_*,V_*)$ is the number of $g$-dimensional subspaces $V'_g$ of $V$ such 
that 

(b) $W_{g-1}\sub V'_g\not\sub W_{n-1}$,
\nl
and

(c) if $g=1\le n-t-1$ then $V'_g\sub V_i,V'_g\not\sub V_{i-1}$ for some $i$ with
$1\le i<j_g$;

(d) if $g\in[2,n-t-1]$ then $V'_g\sub V_i,V'_g\not\sub V_{i-1}$ for some $i$ with 
$j_{g-1}<i<j_g$;

(e) if $g=n-t\ge 2$ then $V'_g\sub V_i,V'_g\not\sub V_{i-1}$ for some $i$ with 
$j_{g-1}<i\le n$.
\nl
Now conditions (c),(d),(e) can be replaced by:

(c${}'$) if $g=1\le n-t-1$ then $V'_g\sub V_{j_g-1}$;

(d${}'$) if $g\in[2,n-t-1]$ then $V'_g\sub V_{j_g-1},V'_g\not\sub V_{j_{g-1}}$;

(e${}'$) if $g=n-t\ge 2$ then $V'_g\not\sub V_{j_{g-1}}$.
\nl
Setting $L=V'_g/W_{g-1}$ we see that $\ph_g(W_*,V_*)$ is the number of lines $L$ in 
$V/W_{g-1}$ such that $L\not\sub W_{n-1}/W_{g-1}$ and

if $g=1\le n-t-1$ then $L\sub V_{j_g-1}/W_{g-1}$;

if $g\in[2,n-t-1]$ then $L\sub V_{j_g-1}/W_{g-1},L\not\sub V_{j_{g-1}}/W_{g-1}$;

if $g=n-t\ge 2$ then $L\not\sub V_{j_{g-1}}/W_{g-1}$.
\nl
Since $W_{n-1}/W_{g-1}$ is a hyperplane in $V/W_{g-1}$, we see that $\ph_g(W_*,V_*)$ is
given by:

$(q^{j_g-g}-q^{j_g-g-1})/(q-1)=q^{j_g-g-1}$ if $g=1\le n-t-1$ and 
$V_{j_g-1}\not\sub W_{n-1}$,

$0$ if $g=1\le n-t-1$ and $V_{j_g-1}\sub W_{n-1}$,

$(q^{j_g-g}-q^{j_g-g-1}-q^{j_{g-1}-g+1}+q^{j_{g-1}-g})/(q-1)=q^{j_g-g-1}-q^{j_{g-1}-g}$
if $g\in[2,n-t-1]$ and $V_{j_{g-1}}\not\sub W_{n-1}$,

$(q^{j_g-g}-q^{j_g-g-1})/(q-1)=q^{j_g-g-1}$ if $g\in[2,n-t-1]$ and 
$V_{j_{g-1}}\sub W_{n-1},V_{j_g-1}\not\sub W_{n-1}$,

$0$ if $g\in[2,n-t-1]$ and $V_{j_g-1}\sub W_{n-1}$,

$(q^{n-g+1}-q^{n-g}-q^{j_{g-1}-g+1}+q^{j_{g-1}-g})/(q-1)=q^{n-g}-q^{j_{g-1}-g}$
if $g=n-t\ge 2$ and $V_{j_{g-1}}\not\sub W_{n-1}$, 

$q^{n-g}$ if $g=n-t\ge 2$ and $V_{j_{g-1}}\sub W_{n-1}$,

$q^{n-g}$ if $g=1=n-t$.
\nl
Now there is a unique $u\in[1,n]$ such that $V_{u-1}\sub W_{n-1},V_u\not\sub W_{n-1}$.
From (a) we see that $u\notin\{j_1,j_2,\do,j_{n-t-1}\}$ (we use that $n-t-1<n-1$). 
Using the formulas above, we can now compute $N=\sum_{g\in[1,n-t]}\ph_g(W_*,V_*)$.

If $u<j_1$ and $n-t\ge 2$ (so that $V_{j_1-1}\not\sub W_{n-1}$) we have

$N=q^{j_1-2}+\sum_{g=2}^{n-t-1}(q^{j_g-g-1}-q^{j_{g-1}-g})+(q^t-q^{j_{n-t-1}-1})=q^t$.
\nl
If $j_{h-1}<u<j_h$ for some $h\in[2,n-t-1]$ (so that 
$V_{j_{h-1}}\sub W_{n-1},V_{j_h-1}\not\sub W_{n-1}$) we have

$N=q^{j_h-h-1}+\sum_{g=h+1}^{n-t-1}(q^{j_g-g-1}-q^{j_{g-1}-g})+(q^t-q^{j_{n-t-1}-1})
=q^t$.
\nl
If $j_{n-t-1}<u$ and $n-t\ge 2$ (so that $V_{j_{n-t-1}}\sub W_{n-1}$) we have $N=q^t$.

If $n-t=1$ we have $N=q^t$.

We see that in any case we have $N=q^t$.

Assume next that $g\in[n-t+1,n]$. If a $V'_g\in Y_g$ exists then there exists
$1\le i_1<i_2<\do<i_{n-t}\le n$ such that

(f) $W_r\sub V_{i_r},W_r\not\sub V_{i_r-1}$ if $r\in[1,n-t]$.
\nl
We see that $\ph_g(W_*,V_*)=0$ if $(W_*,V_*)\notin X_t$. We now assume that 
$(W_*,V_*)\in X_t$ and that $1\le i_1<i_2<\do<i_{n-t}\le n$ is such that (f) holds.
Then $\ph_g(W_*,V_*)$ is the number of $g$-dimensional subspaces $V'_g$ of $V$ such 
that 

$W_{g-1}\sub V'_g\not\sub W_{n-1}$
\nl
that is, the number of lines $L$ in $V/W_{g-1}$ such that $L\not\sub W_{n-1}/W_{g-1}$. 
We see that $\ph_g(W_*,V_*)=q^{n-g}$. Hence
$\sum_{g\in[n-t+1,n]}\ph_g(W_*,V_*)=1+q+q^2+\do+q^{t-1}$.

Summarizing, we see that for $(W_*,V_*)\in\cf\T\cf$, 
$f(W_*,V_*)=\sum_{g=1}^n\ph_g(W_*,V_*)$ is equal to

$1+q+q^2+\do+q^t$ if $(W_*,V_*)\in X_t$,

$q^t$ if $(W_*,V_*)\in X_{t+1}-X_t$,

$0$ if $(W_*,V_*)\notin X_{t+1}$.
\nl
The lemma follows immediately.

\subhead 4\endsubhead
We show that

(a) $q^{1+2+\do+(t-1)}f_t=f_1*(f_1-1)*(f_1-1-q)*\do*(f_1-1-q-q^2-\do-q^{t-2})$
\nl
for $t\in[1,n-1]$ by induction on $t$. For $t=1$ this is clear. Assume that 
$t\in[2,n-1]$ and that (a) holds when $t$ is replaced by $t-1$. Using Lemma 3 we have
$q^{t-1}f_t=(f_1-1-q-q^2-\do-q^{t-2})*f_{t-1}$. Using this and the induction hypothesis
we have
$$\align&q^{1+2+\do+(t-1)}f_t=(f_1-1-q-q^2-\do-q^{t-2})*f_1*(f_1-1)*\\&
*(f_1-1-q)*\do*(f_1-1-q-q^2-\do-q^{t-3}).\endalign$$
This proves (a).

Next we note that $X_{n-1}$ is the set of all $(V'_*,V_*)\in\cf\T\cf$ such that for 
some $i\in[1,n]$ we have $V'_1\sub V_i,V'_1\not\sub V_{i-1}$. 
Thus, $X_{n-1}=\cf\T\cf=X_n$ so that $f_{n-1}=f_n$. Using this and Lemma 3 we see that
$f_1*f_{n-1}=(1+q+q^2+\do+q^{n-1})f_{n-1}$ that is 
$(f_1-1-q-q^2-\do-q^{n-1})f_{n-1}=0$. Hence multiplying both sides of (a) (for $t=n-1$)
by $(f_1-1-q-q^2-\do-q^{n-1})$ we obtain
$$\align&f_1*(f_1-1)*(f_1-1-q)*\do*\\&
*(f_1-1-q-q^2-\do-q^{n-3})*(f_1-1-q-q^2-\do-q^{n-1})=0.\endalign$$
Thus an identity like that in Proposition 2 holds in $\ch$ instead of $H$ (with $f_1,q$
instead of $\t,\qq$). It is known that the algebra $\ch$ may be identified with 
$\CC\ot_{\ZZ[\qq]}H$ (where $\CC$ is regarded as a $\ZZ[\qq]$-algebra via the 
specialization $\qq\m q$) in such a way that $1\ot\t$ is identified with $f_1$. Since 
$q$ can take infinitely many values, the identity in Proposition 2 follows.

\subhead 5\endsubhead
Setting $f_t=0$ for $t>n$, we see that the identity in Lemma 3 remains valid for any 
$t\ge 0$. We see that subspace of $\ch$ spanned by $\{f_t;t\ge 0\}$ coincides with the
subspace spanned by $\{f_1^t;t\ge 0\}$; in particular it is a commutative subring.

\subhead 6\endsubhead
Consider the endomorphism of $\QQ(\qq)\ot_{\ZZ[\qq]}H$ given by left multiplication by 
$\t$. Proposition 2 shows that the eigenvalues of this endomorphism are in 

$\{0,1,1+\qq,1+\qq+\qq^2,\do,1+\qq+\do+\qq^{n-3},1+\qq+\do+\qq^{n-1}\}$. 
\nl
The multiplicity of the eigenvalue $1+\qq+\do+\qq^{k-1}$ is preserved by the
specialization $\qq=1$ hence it is the same as the multiplicity of the eigenvalue $k$ 
for the left multiplication by $\tt$ on $\CC[S_n]$, which by \cite{\DFP} is the number 
of permutations of $[1,n]$ with exactly $k$ fixed points.


\Refs
\widestnumber\key{\DFP}
\ref\key{\DFP}\by P.Diaconis, J.A.Fill and J.Pitman\paper Analysis of top to random
shuffles\jour Combinatorics, probability and computing\vol 1(1992)\pages 135-155\endref
\ref\key{\GW}\by A.M.Garsia and N.Wallach\paper Qsym over Sym is free\yr 2002\endref
\ref\key{\PH}\by R.M.Phatarfod\paper On the matrix occuring in a linear search problem
\jour Jour.Appl.Prob.\vol 28\yr 1991\pages 336-346\endref
\ref\key{\WA}\by N.Wallach\paper Lie algebra cohomology and holomorphic continuation of
generalized Jacquet integrals\jour Advanced Studies in Pure Math.\vol 14(1988)\pages 
123-151\endref
\endRefs
\enddocument